\documentclass[a4paper]{article}
\usepackage{url}
\usepackage{amsmath,amssymb,amsfonts}
\usepackage{algorithm} 
\usepackage{algpseudocode} 
\usepackage{booktabs}
\usepackage{tabularx}
\usepackage{makecell}
\usepackage{placeins}
\usepackage{float}
\usepackage{multirow,subcaption}
\usepackage{tikz}
\usetikzlibrary{shapes.geometric, arrows.meta, positioning}
\usepackage{graphicx,color}
\usepackage{textcomp}

\usepackage{hyperref}
\hypersetup{unicode,
    colorlinks=true}
\usepackage[margin=1in]{geometry} 
\title{Nonlinear Model Order Reduction of Power Grid Networks using Quadratic Manifolds}

\author{Farhana Farooq\thanks{Department of Electrical Engineering, Institute of Technology, University of Kashmir, Srinagar, Jammu \& Kashmir, India 190006 (\url{farhanafarooq19000@gmail.com},\url{danishrafiq@uok.edu.in})} \and Danish Rafiq\footnotemark[1]}

\begin{document}
\maketitle
\begin{abstract}
The increasing size and complexity of modern power systems have led to a high-dimensional mathematical model for transient stability studies, rendering full-scale simulations computationally burdensome.  While dimensionality reduction is essential for reducing this complexity, conventional approaches in power systems predominantly rely on linear projection methods. Such linear subspaces have limited capability for representing the inherently nonlinear swing dynamics of synchronous machines, often resulting in poor approximations and inefficient compression. To address these limitations, this paper introduces a quadratic manifold-based model order reduction (MOR) framework to accelerate the transient dynamic simulations in power systems. The proposed method combines the linear proper orthogonal decomposition (POD) basis with a learned quadratic correction term that minimizes the reconstruction error. This yields a scalable MOR strategy capable of handling strongly nonlinear behaviors, particularly those arising during fast-acting faults, where linear techniques typically fail. The method is tested on a range of benchmark power system models of increasing size and complexity. In addition, we provide a detailed numerical algorithm for constructing the quadratic manifold, along with the corresponding implementation code.    
\end{abstract}
\section{Introduction}
A system is considered large-scale when it exhibits high state dimensionality and is designed to process complex computational tasks and large data volumes \cite{antoulas2005approximation}. Power systems inherently fall into this category due to their numerous interconnected subsystems. As power networks continue to expand, their structural and dynamic complexity increases, rendering analysis and control increasingly demanding in terms of computation \cite{kundur2007power}. The size of a power system is commonly characterized by its number of buses and generators \cite{hamid2023deep}. Contemporary large-scale power systems may include thousands of such components, resulting in mathematical models that are both computationally intensive and challenging to simulate \cite{antoulas2005approximation}. Power system dynamics are inherently nonlinear, further complicating modeling and simulation \cite{nugroho2019characterizing,karlsson2002modelling}. Moreover, storing and processing large volumes of simulation data imposes significant memory and computational burdens, rendering full-order simulations impractical for tasks such as optimization, control design, and probabilistic analysis. Consequently, reduced-order models are essential for approximating the behavior of the full system while preserving its key dynamic characteristics. In power system studies, model reduction is often referred to as \textit{dynamic equivalencing} \cite{chow2013power}. A common approach is to partition the network into a \textit{study} area and an \textit{external} area. The study area contains the components of primary interest and is modeled in full detail, while the external area is simplified—often using linear reduced-order representations—so that its influence on the study area is retained without incurring high computational cost \cite{chaniotis2005model}.
\par
Early model order reduction (MOR) techniques for power systems primarily relied on the network's physical characteristics. One of the earliest approaches was the coherency-based methods \cite{tyuryukanov2020slow, undrill2007construction, germond2007dynamic}, which identify groups of generators that respond similarly to system disturbances and replace each group with an equivalent aggregated model. Although effective when coherent groups are clearly distinguishable, its applicability is limited in systems where generator behavior is less uniform. Subsequent developments introduced synchrony-based methods \cite{ramaswamy2002synchrony, motter2013spontaneous}, singular perturbation techniques \cite{winkleman1981analysis}, and modal analysis approaches \cite{martins2003computing, martins1996computing}, offering improved accuracy but still relying on structural or physical assumptions. As power systems evolved—mainly driven by the increased integration of renewable energy sources—the limitations of the physical interpretation–based reduction became more pronounced. This motivated a shift towards mathematically grounded MOR techniques \cite{moore2003principal, meyer2002fractional, byrnes1995complete}. Methods such as Krylov subspace projection \cite{chaniotis2005model}, balanced truncation \cite{ramirez2015application, benner2013improved}, and extended balanced truncation \cite{sturk2014coherency} have since been successfully applied to large-scale systems. Additional strategies, including measurement-based reduction \cite{chakrabortty2010measurement, meier2003approximation}, ANN-based boundary matching \cite{ma2012hybrid}, heuristic optimization \cite{cepeda2012identification}, and independent component analysis (ICA) \cite{ariff2012coherency}, have further expanded the MOR toolbox.
\par
Despite these advances, most existing MOR methods for power systems still rely on linearization of the external area, which limits their ability to represent the inherently nonlinear dynamics of the high-fidelity model accurately. Nonlinear MOR frameworks have been introduced to mitigate this issue \cite{zhao2014nonlinear, malik2016reduced, osipov2018adaptive}, and parameter-preserving techniques have also been explored \cite{acle2019parameter, scarciotti2016low}. However, these approaches typically require large sets of snapshot data generated under diverse excitation conditions or involve solving computationally demanding \textit{Lyapunov}-type equations. The computational burden becomes even more prohibitive when multiple simulations are needed for varying inputs or parametric changes, restricting their practical applicability to large-scale systems.
\par To overcome these limitations, this paper adopts a nonlinear projection-based MOR strategy, an approach that has received limited attention in power system studies. The key novelty lies in using a quadratic manifold-based reduction framework to capture rotor swing dynamics without requiring repeated full-order simulations. A quadratic mapping operator is identified via a least-squares procedure applied to the residual remaining after projection onto a low-dimensional linear subspace of the full-order model. This yields a compact yet highly accurate surrogate that can replace the full-order model in computationally intensive transient stability studies. The resulting reduced-order formulation enables efficient and reliable simulation of large-scale power system dynamics under both unfaulted and faulted conditions, offering a significant improvement over traditional linear subspace–based reduction techniques. 
\par The remainder of this paper is organized as follows: In section (\ref{sec:PS model}), we describe the dynamic model of power systems. In section (\ref{sec:dim red}), we present the general framework of dimensionality reduction using projective MOR based on linear and nonlinear projections. In section (\ref{sec-results}), we present numerical simulation results for the unfaulted and faulted scenarios. Finally, Section (\ref{sec:con}) concludes the paper. 
\section{Power System Dynamic Model}\label{sec:PS model}
To illustrate the dynamic behavior of large-scale power systems and motivate the need for MOR, we consider the classical nonlinear swing equation, which governs the electromechanical dynamics of synchronous generators \cite{grainger1999power}. This model captures the coupling between rotor angles and electrical power. Under steady-state conditions, the rotor axis of a synchronous machine maintains a fixed position relative to the stator’s rotating magnetic field \cite{fitzgerald2003electric}. The angular displacement between these axes, known as the load or torque angle $\delta$, depends on the machine loading \cite{wadhwa2006electrical}, increasing as the load increases. When a load is suddenly added or removed, the rotor accelerates or decelerates relative to the synchronously rotating stator field, producing electromechanical oscillations. For a system of $n$ coupled oscillators (generators or loads), the swing dynamics at the $i$th node are described by the following nonlinear second-order equation:
\begin{equation}\label{eqn-swing}
\frac{2 \mathbf{J}_i}{\boldsymbol{\omega}_R} \boldsymbol{\ddot{\delta}}_i + \frac{\mathbf{D}_i}{\boldsymbol{\omega}_R} \boldsymbol{\dot{\delta}}_i = \mathbf{F}_i + \mathbf{f}_i(\boldsymbol{\delta)}, \quad \text{for } i=1,\ldots,n,
\end{equation}
where $\boldsymbol{\delta}_i$ denotes the rotor angular position of the $i$th oscillator, ${\omega}_R$ denotes the system's frequency, $\mathbf{J}_i$ and $\mathbf{D}_i$ correspond to the inertia and damping coefficients of the $i$th generator, respectively, $\mathbf{F}_i$ is a constant term associated to the mechanical power and $\mathbf{f}_i$ is the nonlinear term corresponding to the electric demand at the $i$th node, described as:
\begin{equation}
\mathbf{f}_i(\boldsymbol{\delta)} = - \sum_{\substack{j=1 \\ j \neq i}}^n \mathbf{K}_{ij} \sin(\boldsymbol{\delta}_i - \boldsymbol{\delta}_j - \boldsymbol{\gamma}_{ij}),
\end{equation}
where $\mathbf{K}_{ij}$ is the coupling coefficient and $\boldsymbol{\gamma}_{ij}$ is the phase shift (or angle) of the admittance between nodes $i$ and $j$.  
Modeling a single oscillator is quite simple, but in large-scale power systems, modeling is complex due to interconnected subsystems. To address this issue, various modeling techniques have been proposed, including structure-preserving models \cite{bergen2007structure} and adaptive frequency models \cite{acebron1998adaptive}. However, the most commonly used models are the effective network (EN) model and the synchronous motor (SM) model \cite{nishikawa2015comparative}. The EN model assumes loads are constant impedance rather than oscillators. In this case, a network reduction approach is used to reduce the number of dynamic nodes from $n$ to $n_g$, where $n_g$ is the number of generators. The EN model is given as:
\begin{equation}
\frac{2\mathbf{J}_i} {\boldsymbol{\omega}_R} \boldsymbol{\ddot{\delta}}_i + \frac{\mathbf{D}_i} {\boldsymbol{\omega}_R} \boldsymbol{\dot{\delta}}_i + 
\sum_{\substack{j=1 \\ j \neq i}}^{n_g} 
\mathbf{K}^{EN}_{ij} \sin(\boldsymbol{\delta}_i - \boldsymbol{\delta}_j - \boldsymbol{\gamma}^{EN}_{ij}) = \mathbf{F}^{EN}_i, 
\quad i = 1, \ldots, n_g
\label{eq:ENmodel}
\end{equation}

\begin{equation}
\mathbf{F}^{EN}_i := \mathbf{P}^{*}_{g,i} - |\mathbf{E}^{*}_i|^2 \mathbf{G}^{EN}_{i,j}, \quad
\mathbf{K}^{EN}_{ij} := \left| \mathbf{E}^{}_i \mathbf{E}^{}_j \mathbf{Y}^{EN}_{i,j} \right|, \quad
\boldsymbol{\gamma}^{EN}_{ij} := \boldsymbol{\lambda}^{EN}_{ij} - \frac{\pi}{2},
\end{equation}
where  $\mathbf{P}^{*}_{g,i}$ denotes the electrical power output of $i$th generator at the balanced operating point, $\mathbf{Y}^{EN}_{i,j}$ represents the corresponding admittance matrix entry between the generator terminal bus, obtained after eliminating the load buses via \textit{kron} reduction \cite{wang2023modelling}, $\mathbf{G}^{EN}_{i,j}$ denotes the real component. At the same time, $\boldsymbol{\lambda}^{EN}_{ij}$ corresponds to its phase angle. Due to the reduction process, the admittance network in the EN model does not retain the original full graph structure. $\mathbf{E}^{*}_i$ represents the generator's internal behind the transient reactance.  On the other hand, the SM modelling approach treats every load as a synchronous motor with negative output power \cite{rafiq2022synergistic}. Unlike the EN model, the SM model fully preserves the network's graph topology \cite{sauer2017power}. The modeling provides the complete dynamic behavior of all main grid components, which are essential for stability analysis. The SM model equations are obtained as:
\begin{equation}
\frac{2\mathbf{J}_i} {\boldsymbol{\omega}_R} \boldsymbol{\ddot{\delta}}_i + \frac{\mathbf{D}_i} {\boldsymbol{\omega}_R} \boldsymbol{\dot{\delta}}_i + 
\sum_{\substack{j=1 \\ j \neq i}}^{n} 
\mathbf{K}^{SM}_{ij} \sin(\boldsymbol{\delta}_i - \boldsymbol{\delta}_j - \boldsymbol{\gamma}^{SM}_{ij}) = \mathbf{F}^{SM}_i, 
\quad i = 1, \ldots, n,
\label{eq:SMmodel}
\end{equation}
\begin{equation}
\mathbf{K}^{SM}_{ij} := \left| \mathbf{E}^{}_i \mathbf{E}^{}_j \mathbf{Y}^{SM}_{i,j} \right|, \quad \mathbf{Y}^{SM}_{i,j} = \left| \mathbf{Y}^{SM}_{i,j} \right| \mathbf{e}^{j \boldsymbol{\lambda}^{SM}_{ij}}, \quad
\boldsymbol{\gamma}^{SM}_{ij} := \boldsymbol{\lambda}^{SM}_{ij} - \frac{\pi}{2}.
\end{equation}
Here $\mathbf{e}^{j\boldsymbol{\lambda}^{SM}_{ij}}$ represents the complex exponential form with $j = \sqrt{-1}$,
expressing the phase angle of the admittance element. 
Although the parameters are expressed using the same equations as in the EN model, the physical interpretation of the symbols changes. In the SM formulation, the quantities represent the full synchronous machine network, whereas in the EN model, they describe the reduced equivalent system.\par 
To facilitate numerical simulation and enable the application of MOR techniques, the second-order swing equations are rewritten in first-order nonlinear state-space form as: 
 \begin{align}\label{eqn-swing1order}
 	\boldsymbol{\dot{\delta}}_i &= \boldsymbol{\omega}_i, \\
 	\boldsymbol{\dot{\omega}}_i &= -\frac{\mathbf{D}_i}{2\mathbf{J}_i}\boldsymbol{\omega}_i + \frac{\boldsymbol{\omega}_R}{2\mathbf{J}_i}\mathbf{F}_i + \frac{\boldsymbol{\omega}_R}{2\mathbf{J}_i}\mathbf{f}_i(\boldsymbol{\delta}),
 \quad \text{for } i = 1, \ldots, n
\end{align}
 where \(\boldsymbol{\omega}_i\) is the frequency deviation of the \(i\)th generator. Transforming the dynamics into this first-order form results in a total of $2n$ state variables for a system with $n$ generators. Consequently, a large-scale power networks, often comprising hundreds or thousands of generators, yield extremely high-dimensional nonlinear dynamical models. Simulating such models requires solving a large set of coupled nonlinear differential equations, which quickly becomes computationally prohibitive, particularly for transient stability studies, parametric analyses, or scenarios that require repeated simulations. In the next section, we discuss the MOR procedure for obtaining low-order models for such systems.

 \section{Dimensionality reduction using projective MOR} \label{sec:dim red}
Dimensionality reduction aims to obtain a low-dimensional representation of a high-dimensional system. This is often achieved using \textit{projection}, i.e., by projecting the full order model dynamics onto a low-dimensional subspace and then describing the reduced dynamics within that subspace \cite{varona2020model}. In the following, we discuss two types of projections.  
\subsection{Linear projection based MOR (Linear Manifold)}
Consider a nonlinear, continuous-time model described in state-space as follows:
\begin{equation}
\begin{split}\label{eq:Fom equation}
		\mathbf{E} \mathbf{\dot{x}}(t)&= \mathbf{f}(\mathbf{x}(t), \mathbf{u}(t)), \quad \mathbf{x}(0) = \mathbf{x}_0, \\
	\mathbf{y}(t) &= \mathbf{h}(\mathbf{x}(t)),
\end{split} 
\end{equation}
where $\mathbf{x}(t) \in \mathbb{R}^{n}$ is the state-vector comprising the state-variables with $\mathbf{x}_0$ as the initial conditions, $\mathbf{u}(t) \in \mathbb{R}^{m}$
is the input vector, $\mathbf{y}(t) \in \mathbb{R}^{p}$ is the output vector, $ \mathbf{E} \in \mathbb{R}^{n \times n}$ is the non-singular descriptor matrix, $\mathbf{f}(\mathbf{x},\mathbf{u}) : \mathbb{R}^{n} \times \mathbb{R}^{m} \to \mathbb{R}^{n}$ and $\mathbf{h}(\mathbf{x}) : \mathbb{R}^{n} \to \mathbb{R}^{p}$ are two nonlinear functions that describe the evolution of the state. We refer to (\ref{eq:Fom equation}) as the full order model (FOM), which represents the true description of the system.\par
An established and successful way to reduce the dimensionality of (\ref{eq:Fom equation}) is to apply the classical \textit{Petrov-Galerkin} projection, i.e., the high-dimensional state vector $\mathbf{x}(t)$ is approximated by a low-dimensional variable $\mathbf{x}_r(t)$ given as:
\begin{equation}
	\mathbf{x}(t) = \mathbf{V} \mathbf{x}_r(t) + e(t)
    \label{eq:approx}
\end{equation}
where  $\mathbf{V} \in \mathbb{R}^{n\times r}$ represents the reduced basis, that we intend to discover, $\mathbf{x}_r(t) \in \mathbb{R}^{r}$, and $e(t) \in \mathbb{R}^{n}$ denotes the approximation error due to projection. 
Using the approximation (\ref{eq:approx}) into the FOM model (\ref{eq:Fom equation})
leads to an overdetermined system given as:
\begin{equation}
\label{eq:Fom_approx}
		\mathbf{E}\mathbf{V\dot{x}_r}(t)- \mathbf{f}(\mathbf{Vx}_r(t), \mathbf{u}(t))-\boldsymbol{\varepsilon}(t)=\mathbf{0}
\end{equation}
where $\boldsymbol{\varepsilon}(t)$ is the residual term. Now, pre-multiplying (\ref{eq:Fom_approx}) by the projector $\boldsymbol{\pi}$ leads to:
\begin{equation}
	\boldsymbol{\pi} \Big( \underbrace{\mathbf{E}\mathbf{V}\dot{\mathbf{ x}}_r(t)-{\mathbf{f\big({V}{x}}_r(t),\mathbf{u}(t)}\big)}_{\boldsymbol{\xi}(\mathbf{Vx}_r(t), \mathbf{u}(t))}-\boldsymbol{\varepsilon}(t)\Big) =  \mathbf{0} \implies \boldsymbol{\pi} \Big({\boldsymbol{\xi}(\mathbf{Vx}_r(t), \mathbf{u}(t))}-\boldsymbol{\varepsilon}(t)\Big)= \mathbf{0}.
	\end{equation} 
	Enforcing the \textit{Petrov-Galerkin} condition, i.e., $\mathbf{W}^{\top}{\boldsymbol{\varepsilon}}(t)={0}$, which implies $\boldsymbol{{\pi}\,{\boldsymbol{\varepsilon}}}(t)={0}$, the residual then becomes zero and only the term ${\boldsymbol{\pi}}\,{\boldsymbol{\xi}}\big(\mathbf{{V}{x}_{\rm r}}(t),\, \mathbf{u}(t)\big)={0}$ remains. This finally yields the reduced order model (ROM) given as:
	\begin{align}
    \mathbf{E}_{r}\,\dot{\mathbf{x}}_r(t) &= \mathbf{{W}^{\top}}\mathbf{f}_r\big(\mathbf{V}{\mathbf{x}}_r(t),\,\mathbf{u}(t)\big),\quad {\mathbf{x}}_r(0)={\mathbf{x}}_{{r}0}, \\
		\mathbf{y}_r(t) &= \mathbf{h}_r\big(\mathbf{V}{\mathbf{x}}_r(t)\big), 
	\end{align}
	where $\mathbf{E}_{r}={\mathbf{W}}^{\top}{\mathbf{E}}{\mathbf{V}}$ is the reduced descriptor matrix, $\mathbf{W} \in \mathbb{R}^{n \times r}$ is another projection matrix, $\mathbf{x}_r(0)=({\mathbf{W}}^{\top}{\mathbf{E}}{\mathbf{V}})^{-1}{\mathbf{W}}^{\top}\mathbf{E}{\mathbf{x}}_0$ and the reduced nonlinear function $\mathbf{f}_r\big(\mathbf{x}_r(t),\mathbf{u}(t)\big)=\mathbf{W}^{\top}\mathbf{f}\big(\mathbf{V}\mathbf{x}_r(t), \mathbf{u}(t)\big)$ with $\mathbf{f}_r(\mathbf{x}_r,\mathbf{u}):\mathbb{R}^{r}\times \mathbb{R}^{m}\to\mathbb{R}^{r}$.
    \par 
The most successful linear projection-based technique is the Proper orthogonal decomposition (POD) \cite{hinze2005proper,moore2003principal}. The POD method aims at finding a basis that optimally approximates the space spanned by an arbitrary set of data points, also called \textit{snapshots} that represent discrete samples $\{\mathbf{x}(t_k)_{k=1}^{k}\}$ of the FOM state trajectory. Given a set of snapshots $\mathbf{X} = [\mathbf{x}(t_1), ..., \mathbf{x}(t_k)]$ where $\mathbf{X} \in \mathbb{R}^{n \times k}$ of the FOM state trajectory, the goal of POD is to find $\mathbf{V}=[\mathbf{v}_1,...,\mathbf{v}_r] \in \mathbb{R}^{n \times r}$ such that:
\begin{equation}
	\min_{V} \sum_{i=1}^{k} \left\| \mathbf{x}(t_i) - \mathbf{V} \mathbf{x}_{r}(t_i) \right\|_2^2 
	= \min_{\{v_j\}} \sum_{i=1}^{k} \left\| \mathbf{x}(t_i) - \sum_{j=1}^{r} \mathbf{v}_j \mathbf{v}_j^{T} \mathbf{x}(t_i) \right\|_2^2, \quad s.t.~~ \mathbf{v}_i^{\top}\mathbf{v}_{j}=\delta_{ij}  
	\label{eq:pod_opt}
\end{equation}
where $\mathbf{x}_{r}(t_i) \in \mathbb{R}^{r}$ is the coefficient vector used to minimize the above expression. The optimization minimizes the reconstruction error by projecting the high-dimensional data onto a reduced basis.  Actually, this is the least squares minimization problem, and the minimization problem can be solved by applying a singular value decomposition (SVD) of the snapshot matrix $\mathbf{X}$ \cite{taira2017modal}. The reduced basis $\mathbf{V}$ is finally constructed by taking the first $r$ left singular vectors corresponding to the largest singular values $\xi_i$. A typical approach for choosing $r$ is that the approximation in the $r$-dimensional basis $\mathbf{V} \in \mathbb{R}^{n \times r}$ yields
\begin{equation}\label{eqn-energy_linear}
    \dfrac{\| \mathbf{VV}^\top (\mathbf{X} - \mathbf{X}_{\mathrm{ref}}) \|_F^2}{\| \mathbf{X} - \mathbf{X}{\mathrm{ref}} \|_F^2}
= \frac{\sum_{i=1}^{r} \sigma_i^2}{\sum_{i=1}^{n} \sigma_i^2} > \kappa,
\end{equation}
where $\kappa$ is a user-specified tolerance, $\mathbf{X}_{\text{ref}}$ is the reference matrix that contains the reference state $\mathbf{x}_{\text{ref}}$ that shifts the original data, and the $\sigma^2$ are the squared singular values of the shifted data matrix. Equation (\ref{eqn-energy_linear}) is also called the \textit{relative cumulative energy} captured by the $r$ leading POD modes.
\par 
POD is a straightforward, data-driven technique for obtaining the columns of the basis $\mathbf{V}$. For this reason, this is the well-known nonlinear reduction technique. However, POD does not effectively capture the system's nonlinear behavior. This is true, especially when significant nonlinearities are present. In this case, a large number of POD modes may be required for the linear projection, increasing computational cost and reducing the effectiveness of the reduction. For systems exhibiting mild nonlinear behavior, the trajectory piece-wise linear (TPWL) method extends this idea by linearizing systems at several points along a representative trajectory \cite{white2003trajectory, rewienski2003trajectory}. Each linearized model captures the local dynamics around a specific operating point, and their combination provides an accurate global approximation within the linear manifold framework. Moreover, for large-scale systems, linear projection often results in higher CPU times during ROM simulations. To address these issues, hyper-reduction techniques such as the discrete empirical interpolation method (DEIM) \cite{chaturantabut2010nonlinear} have been proposed. The DEIM method improves the POD approximation by reducing the complexity of the nonlinear term to a level proportional to the number of reduced variables \cite{chaturantabut2009discrete}, thereby enhancing computational efficiency. However, even with DEIM, the fundamental limitation remains that POD relies on a linear subspace, which restricts its ability to capture nonlinear dynamics effectively. This limitation motivates the use of nonlinear projection-based techniques, which are discussed next.  
\subsection{Nonlinear Projection based MOR (Quadratic Manifold)}
\begin{figure}[H]
\includegraphics[width=\linewidth]{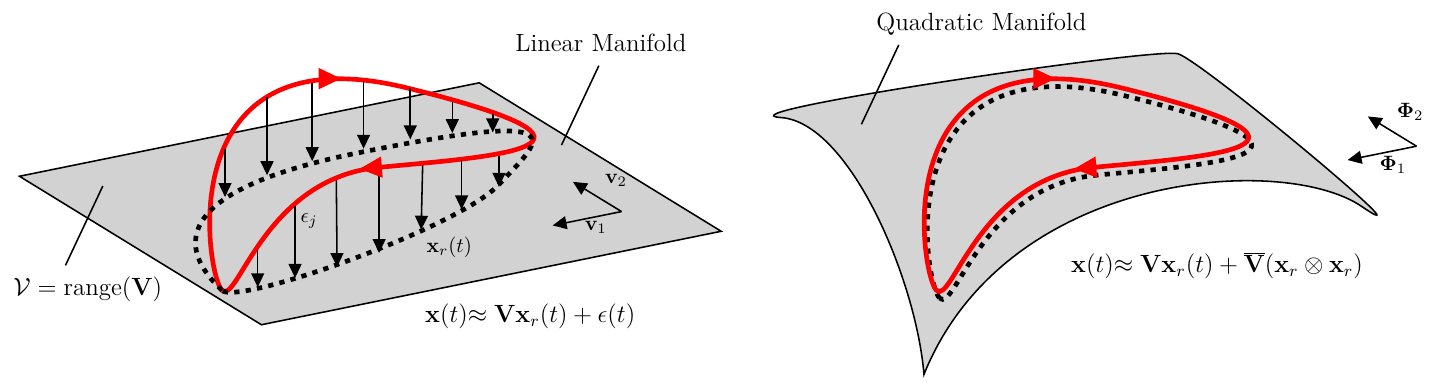}
\caption{Graphical illustrations of linear and nonlinear projection-based model order reduction techniques. The solid red line represents the true solution trajectory, and the dashed black line represents the approximated solution captured on the manifolds.}
\label{fig-manifold}
\end{figure}

Now consider the following approximation of the state vector $\mathbf{x}(t)$ given as:
\begin{equation}
	\mathbf{x}(t) \approx \mathbf{x}_{\mathrm{ref}} + \mathbf{\Phi} (\mathbf{x}_r(t)),
\end{equation}
where $\mathbf{\Phi}:\mathbb{R}^{r} \mapsto \mathbb{R}^{n}$ is a nonlinear function. Notice that there are many possible choices for $\mathbf{\Phi}$, namely cubic, septical, and sinusoidal, but in this work, we restrict the use of the quadratic function. Using a nonlinear projection basis means projecting the system dynamics onto a manifold rather than a linear subspace. This has the advantage that a trial manifold can minimize the residual remaining after projection onto a low-dimensional linear subspace. This is illustrated in Fig. (\ref{fig-manifold}). 
In the figure, the linear manifold (illustrated on the left) is represented by the system states as a linear combination of basis vectors $\mathbf{v}_1$ and $\mathbf{v}_2$. The projection takes place onto a flat subspace, denoted as $\mathcal{V} = range(\mathbf{V})$, where the state vector $\mathbf{x}(t)$ is approximated, as per (\ref{eq:approx}). This method efficiently captures the system's dominant linear dynamics but struggles to accurately describe strongly nonlinear effects. In contrast, the quadratic or nonlinear manifold approach (illustrated on the right) improves the approximation by using a curved subspace rather than a flat one. In this case, the high-dimensional state $\mathbf{x}(t)$ is expressed as:   
\begin{equation}
	\mathbf{x}(t) \approx \boldsymbol{\Phi}(\mathbf{x}_r(t)) :=
	\mathbf{x}_{\mathrm{ref}} + \mathbf{V} \mathbf{x}_r(t) + \mathbf{\bar{V}} \, (\mathbf{x}_r(t) \otimes
	\mathbf{x}_r(t)),
	\label{eq:quad_manifold}
\end{equation}
where $\mathbf{\bar{V}} \in \mathbb{R}^{n \times r^{2}}$ is a data-driven
quadratic mapping operator and $\otimes$ denotes the \textit{Kronecker}
product\footnote{The operator $\otimes$ denotes the Kronecker product which for a column vector $x = [x_1, x_2, ..., x_n]^\top$ is given by $x \otimes x = [x_1^2, x_1x_2, ..., x_1x_n,x_2x_1,x_2^2,..,x_n^2]^\top \in \mathbb{R}^{n^2}$ }. The use of an additional quadratic term enriches the approximation
space with directions in the orthogonal complement of $\mathbf{V}$, thereby
capturing part of the projection error without introducing extra
degrees of freedom in the reduced coordinates $\mathbf{x}_r(t)$. 
The operator $\mathbf{\bar{V}}$ is inferred from snapshot data by solving the
least-squares problem:
\begin{equation}
	\mathbf{\bar{V}} = \arg\min_{\bar{V} \in \mathbb{R}^{n \times r^2}}
	 \sum_{j=1}^k
	\left\|
	\mathbf{x}^j - \mathbf{x}_{\mathrm{ref}} - \mathbf{V} \mathbf{x}_{r}^j
	- \mathbf{\bar{V}}(\mathbf{x}_{r}^j \otimes \mathbf{x}_{r}^j)
	\right\|_2^2.
    \label{eqn-opt}
\end{equation}
This implies that the constructed mapping operator fulfills:
\begin{equation}
    \mathbf{x}_r\approx \mathbf{x}_{\mathrm{ref}}+\mathbf{V}(\mathbf{x}_{r}^j) +\mathbf{\bar{V}}(\mathbf{x}_{r}^j \otimes \mathbf{x}_{r}^j),\quad j=1,...,k, 
\end{equation}
at the optimum of the objective of (\ref{eqn-opt}). Here, $\mathbf{V}$ is obtained via POD and $\mathbf{x}_{r}^j=\mathbf{V}^T(\mathbf{x}^j-\mathbf{x}_{\mathrm{ref}})$ is the $j$th snapshot of reduced variable $\mathbf{x}_{r}^j$. The above optimization problem can also be reformulated in the Frobenius norm as:
\begin{equation}
\mathbf	{\bar{V}} = \arg\min_{\bar{V} \in \mathbb{R}^{n \times r^2}}
	\frac{1}{2} \big\| \mathbf{P}^\top \mathbf{\bar{V}^\top} - \boldsymbol{\mathcal{E}}^\top \big\|_F^2,
    \label{eqn-opt2}
\end{equation}
where $\mathcal{E} = (\mathbf{I} - \mathbf{VV}^\top)(\mathbf{x}- \mathbf{x}_{\mathrm{ref}}):=\mathbf{V}_{\bot}\mathbf{V}_{\bot}^\top(\mathbf{X-\mathbf{X}_{\mathrm{ref}}}) \in \mathbb{R}^{n\times k}$ is the linear
projection error, and
\begin{equation}
	\mathbf{P := \begin{bmatrix}
	|&|& &|\\
		\mathbf{x}_{r}^1 \otimes \mathbf{x}_{r}^1 &
		\mathbf{x}_{r}^2 \otimes \mathbf{x}_{r}^2 & 
		\mathbf{\cdots} &
		\mathbf{x}_{r}^k \otimes \mathbf{x}_{r}^k\\|&|& &|
	\end{bmatrix}} \in \mathbf{R}^{r^2 \times k}.
\end{equation}
Here, the quadratic basis matrix $\mathbf{\bar{V}}$ is obtained by solving the optimization problem (\ref{eqn-opt2}), which reduces the error due to the linear-subspace approximation. The equation (\ref{eqn-opt2}) can also be represented as:
\begin{equation}
\mathbf{P}^{\top}\mathbf{\overline{V}}^{\top} - \boldsymbol{\mathcal{E}}^{\top} :=
	\underbrace{
	\begin{pmatrix}
		(\mathbf{x}_{r}^1 \otimes \mathbf{x}_{r}^1 )^{\top} \\
		(\mathbf{x}_{r}^2  \otimes \mathbf{x}_{r}^2 )^{\top} \\
		\vdots \\
		(\mathbf{x}_{r}^k  \otimes \mathbf{x}_{r}^k )^{\top}
	\end{pmatrix}
}_{k \times r^2}
\underbrace{
	\begin{pmatrix}
\mathbf{\overline{v}}_1^{\top} \\
		\mathbf{\overline{v}}_2^{\top} \\
		\vdots \\
		\mathbf{\overline{v}}_{r^{2}}^{\top}
	\end{pmatrix}
}_{r^2 \times n}
	\underbrace{
	-\begin{pmatrix}
		\boldsymbol{\varepsilon}_1^{\top} \\
		\boldsymbol{\varepsilon}_2^{\top} \\
		\vdots \\
		\boldsymbol{\varepsilon}_n^{\top}
	\end{pmatrix}
}_{n \times n},
	\label{eq:12}
\end{equation}
which is an overdetermined linear least-squares problem  whenever ${n} > {r}^2$. By removing redundant terms in the Kronecker products $(\mathbf{x}_{r}^j \otimes \mathbf{x}_{r}^j)$, the number 
of columns in $\mathbf{\bar{V}}$ can be reduced from ${r}^2$ to 
${r(r+1)}/{2}$. Consequently, the condition for overdetermindedness 
changes to $k > {r(r+1)}/{2}$. It is worth mentioning here that, in 
practice, we often encounter this situation, since the number of available snapshots $k$ is usually much larger 
than the reduced dimension $r$. After reducing the dimension of $\mathbf{\bar{V}}$ from $n \times r^2$ to $n \times r(r+1)/2$, (\ref{eqn-opt2}) can be solved through normal equation given as:
\begin{equation}
\mathbf{\bar{V}}^{\top} =(\mathbf{PP}^{\top})^{-1}\mathbf{P}\mathcal{E}^{\top} \implies \mathbf{\bar{V}}=\mathcal{E}\mathbf{P}^{\top}(\mathbf{PP}^\top)^{-1} \in \mathbb{R}^{n \times r(r+1)/2}.
\end{equation}
Using the value of  $\boldsymbol{\mathcal{E}}$, we have:
\begin{equation}
    \mathbf{\bar{V}}=(\mathbf{I}-\mathbf{VV}^\top)(\mathbf{X}-\mathbf{X}_{\mathrm{ref}})\mathbf{P}^{\top}(\mathbf{PP}^{\top})^{-1}=\mathbf{V}_{\bot}[\mathbf{V}_{\bot}^{\top}\mathbf{X}-\mathbf{X}_{\mathrm{ref}})\mathbf{P}^{\top}(\mathbf{PP}^{\top})^{-1}] \in \mathbb{R}^{n \times r(r+1)/2},
\end{equation}
which implies that each column of $\mathbf{\bar{V}}$ is in the column space of $\mathbf{V}_{\bot}$, such that $\mathbf{V}^{\top}\mathbf{\bar{V}}=\mathbf{0}$ holds. Finally, to avoid noise amplification with relatively little training data, least-squares regularization is adopted to promote solutions that fit the data well. A commonly used practice is to use a Frobenius regularization, which yields:
\begin{equation}
	\mathbf{\bar{V}} := \arg\min_{\bar{V} \in \mathbb{R}^{n \times r(r+1)/2}}
	\left(
	\dfrac{1}{2} \big\| \mathbf{P}^\top \mathbf{\bar{V}}^\top - \boldsymbol{\mathcal{E}}^\top \big\|_F^2
	+ \dfrac{\lambda}{2} \| \mathbf{\bar{V}} \|_F^2
	\right),
	\label{eq:eq15}
\end{equation}
where ${\mathbf{\lambda}} > 0$ is a regularization parameter. The explicit solution then becomes:
\begin{equation}
	\mathbf{\bar{V}} = \mathcal{E} \mathbf{P}^\top (\mathbf{PP}^\top + \lambda \mathbf{I})^{-1} \in \mathbb{R}^{n \times r(r+1)/2}.
\end{equation}
The overall steps for obtaining the matrices $\mathbf{V}$ and $\mathbf{\bar{V}}$ are summarized in Algorithm (\ref{alg:quadratic_manifold}). Similar to (\ref{eqn-energy_linear}), we can define the ``retained energy" metric for the quadratic manifold case, for choosing the reduced basis dimension as:
\begin{equation}\label{eqn-energy-nonlinear}
    \frac{\| \mathbf{VV}^\top (\mathbf{X} - \mathbf{X}_{\mathrm{ref}}) 
+ \mathbf{\overline{V}}\,(\mathbf{V}^\top (\mathbf{X} - \mathbf{X}_{\mathrm{ref}}) \odot \mathbf{V}^\top (\mathbf{X} - \mathbf{X}_{\mathrm{ref}})) \|_F^2}
{\| \mathbf{X} - \mathbf{X}_{\mathrm{ref}} \|_F^2} > \kappa,
\end{equation}
where $\otimes$ is the column-wise Kronecker product of two matrices.
\begin{algorithm}[H]
	\caption{Nonlinear Projection based on Quadratic Manifolds}
	\label{alg:quadratic_manifold}
     \textbf{Input:} Snapshot matrix $\mathbf{X} \in \mathbb{R}^{n \times r}$.\\
     \textbf{Result:} Projection basis $\mathbf{V}, \mathbf{\bar{V}}$
	\begin{algorithmic}[1]
		\item[] \textit{\{Obtain the linear basis using POD\}}
		\State Shift the snapshot data relative to the initial reference state
		\State Perform the Singular Value Decomposition of $\mathbf{X - X_{\text{ref}}}$
		\State $r \leftarrow$ Select reduced basis dimension (\ref{eqn-energy_linear})        
		\State $\mathbf{V} \leftarrow$ Retain the $r$ principal left singular vectors of $\mathbf{X - X_{\text{ref}}}$
		\item[] \textit{\{Obtain the quadratic manifold basis\}}
		\State Express the states in reduced basis: $\mathbf{X}_r \leftarrow \mathbf{V}^{\top}(\mathbf{X} - \mathbf{X}_{\text{ref}})$
		\State $\boldsymbol{\mathcal{E}}, \mathbf{P} \leftarrow$ Compute residual error and quadratic data matrix
		
		\State $\mathbf{\bar{V}} \leftarrow$ Solve regularized least-squares problem using (\ref{eq:eq15})
		\State $\boldsymbol{\lambda} \leftarrow$ Set the regularization parameter
	\end{algorithmic}
\end{algorithm}
The quadratic manifold procedure can be regarded as a
data-driven \emph{closure modeling} technique, in which the effect of the discarded modes is incorporated through the quadratic correction
term $\mathbf{P}(\mathbf{x}_r\otimes \mathbf{x}_r)$. Unlike the conventional linear projection, this methodology establishes a quadratic relationship between the reduced and full-order coordinates, thereby enabling the reduced-order model to account for higher-order nonlinear effects without increasing the reduced dimension. This concept is related to the nonlinear manifold-based reduction strategy proposed in the \textit{maniMOR} \cite{gu2011model}, which represents system trajectories on a nonlinear manifold rather than a purely linear subspace. Importantly, while the manifold
introduces additional directions into the approximation, the
dimension of the reduced coordinates remains $r$, thus avoiding
increased complexity in the reduced-order system. In the next section, we will demonstrate the application of the quadratic manifold for power grid models.

\section{Numerical Simulations}\label{sec-results}
The simulations were performed in MATLAB R2024a on a workstation equipped with an 11th Gen Intel\textsuperscript{\textregistered} Core\texttrademark{} i7-14700K (3.40 GHz) processor and 32 GB RAM. Generator, bus, and branch data were imported using the MATPOWER library \cite{zimmerman2016matpower} to construct the matrices $\mathbf{J}, \mathbf{D}, \mathbf{F}$ and the nonlinear function $\mathbf{f}_i$ for the swing equation (\ref{eqn-swing}). The EN modeling approach described in Section (\ref{sec:PS model}) was adopted. The swing dynamics were integrated using the implicit \textit{Euler's} method \cite{quarteroni2006numerical} with a fixed step size of $\Delta t = 0.01s$. To evaluate the accuracy and computational efficiency of the quadratic manifold reduction, three benchmark IEEE test systems were considered: the IEEE 118-bus, IEEE 300-bus, and Polish 2736-bus models. Performance was assessed relative to the linear POD-based manifold under both normal operating conditions and abnormal (faulted) scenarios to examine the overall robustness of the reduction methods.
\subsection{Case I: Systems under normal conditions}
\begin{figure*}[t] 
	\centering
	\includegraphics[width=\textwidth]{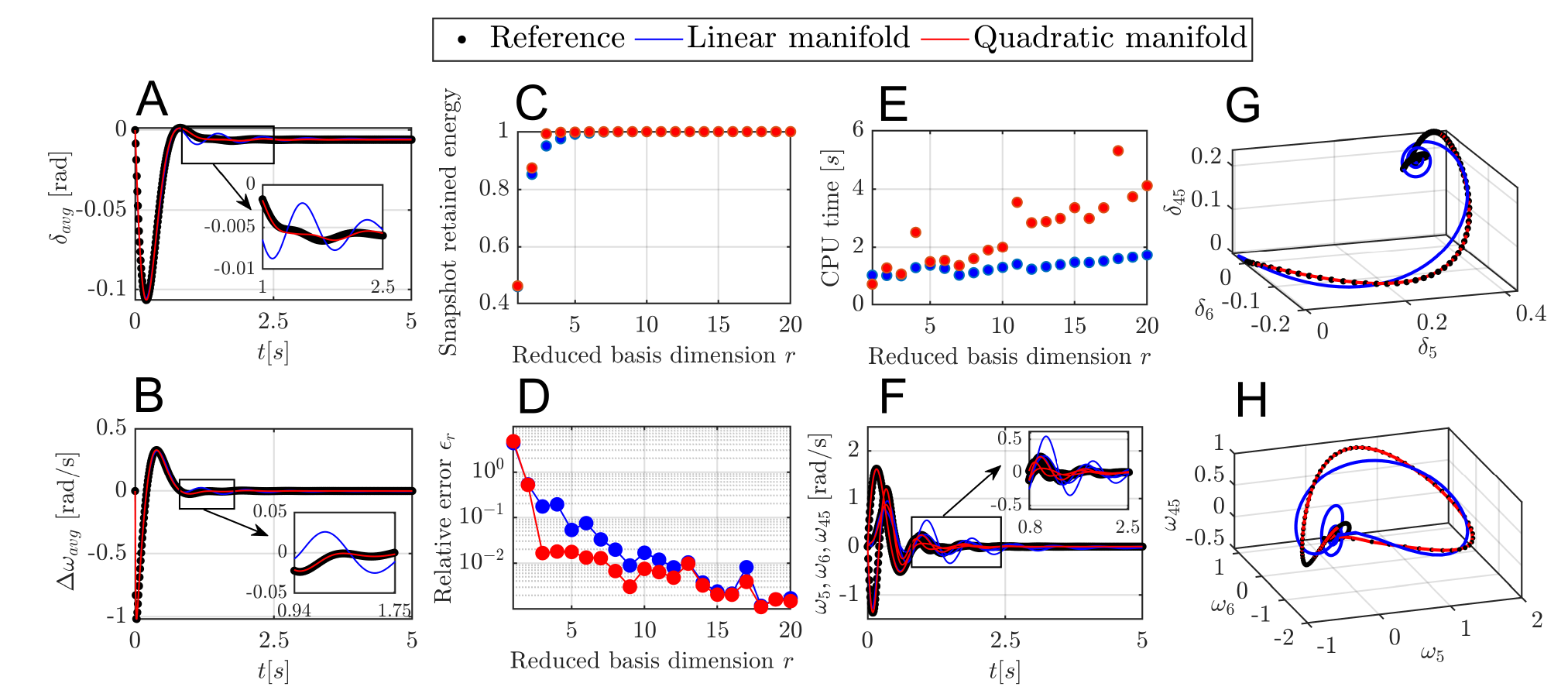} 
	\caption{Simulation results of IEEE 118 bus system: (A) variation of $\delta_{avg}$ with respect to time; (B) variation of $\Delta\omega_{avg}$ with respect to time; (C) snapshot energy captured by the system at various reduced basis dimensions; (D) comparison of relative error for linear and quadratic manifold ROMs; (E) comparison of the CPU times for increasing basis dimension; (F) comparison of the output response for three individual generators with respect to time; (G) comparison os system trajectory onto a linear manifold and quadratic manifold; (H) evolution of angular velocities with respect to time}
	\label{fig:118}
\end{figure*}
\begin{figure*}[t] 
	\centering
	\includegraphics[width=\textwidth]{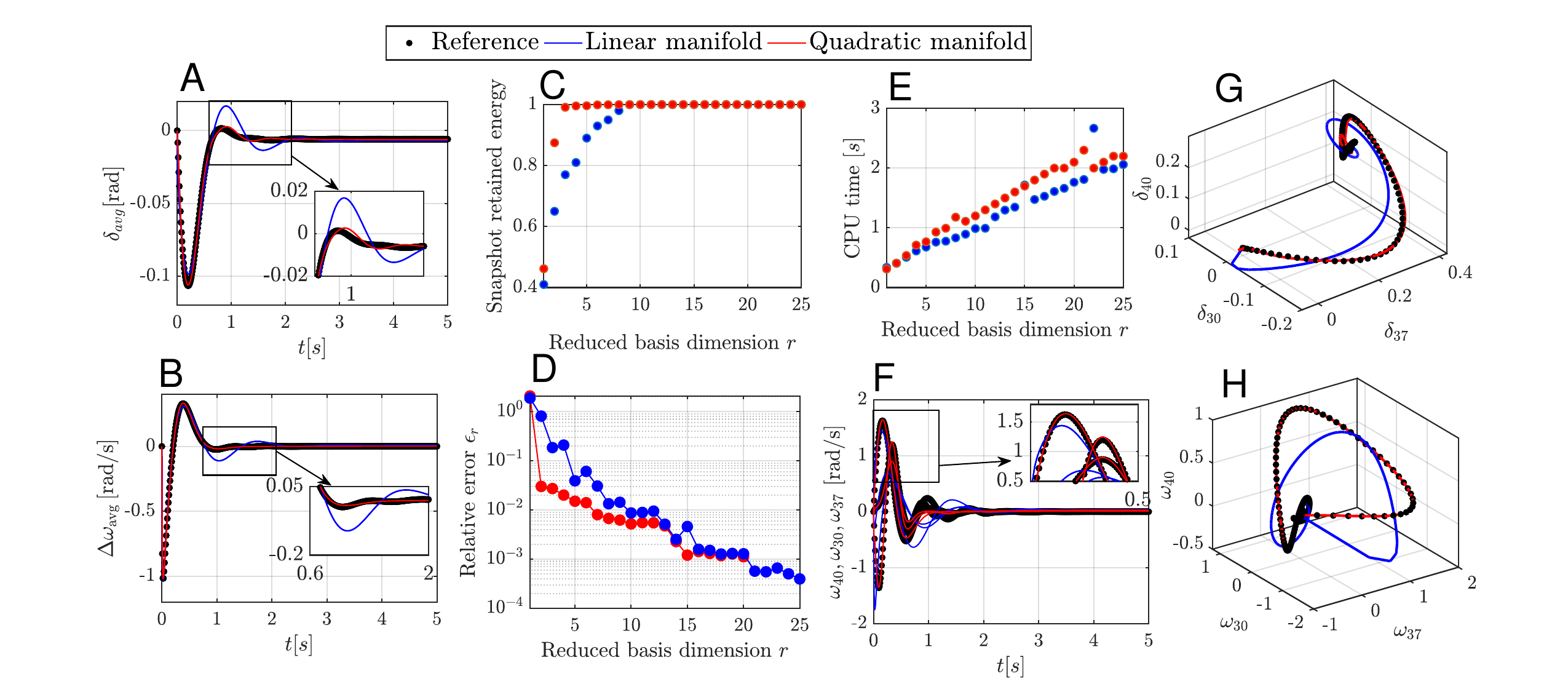} 
	\caption{Simulation results of IEEE 300 bus system: (A) variation of $\delta_{avg}$ with respect to time; (B) variation of $\Delta\omega_{avg}$ with respect to time; (C) snapshot energy captured by the system at various reduced basis dimensions; (D) comparison of relative error for linear and quadratic manifold ROMs; (E) comparison of the CPU times for increasing basis dimension; (F) comparison of the output response for three individual generators with respect to time; (G) comparison os system trajectory onto a linear manifold and quadratic manifold; (H) evolution of angular velocities with respect to time}
	\label{fig:300}
\end{figure*}
\begin{figure*}[t] 
	\centering
	\includegraphics[width=\textwidth]{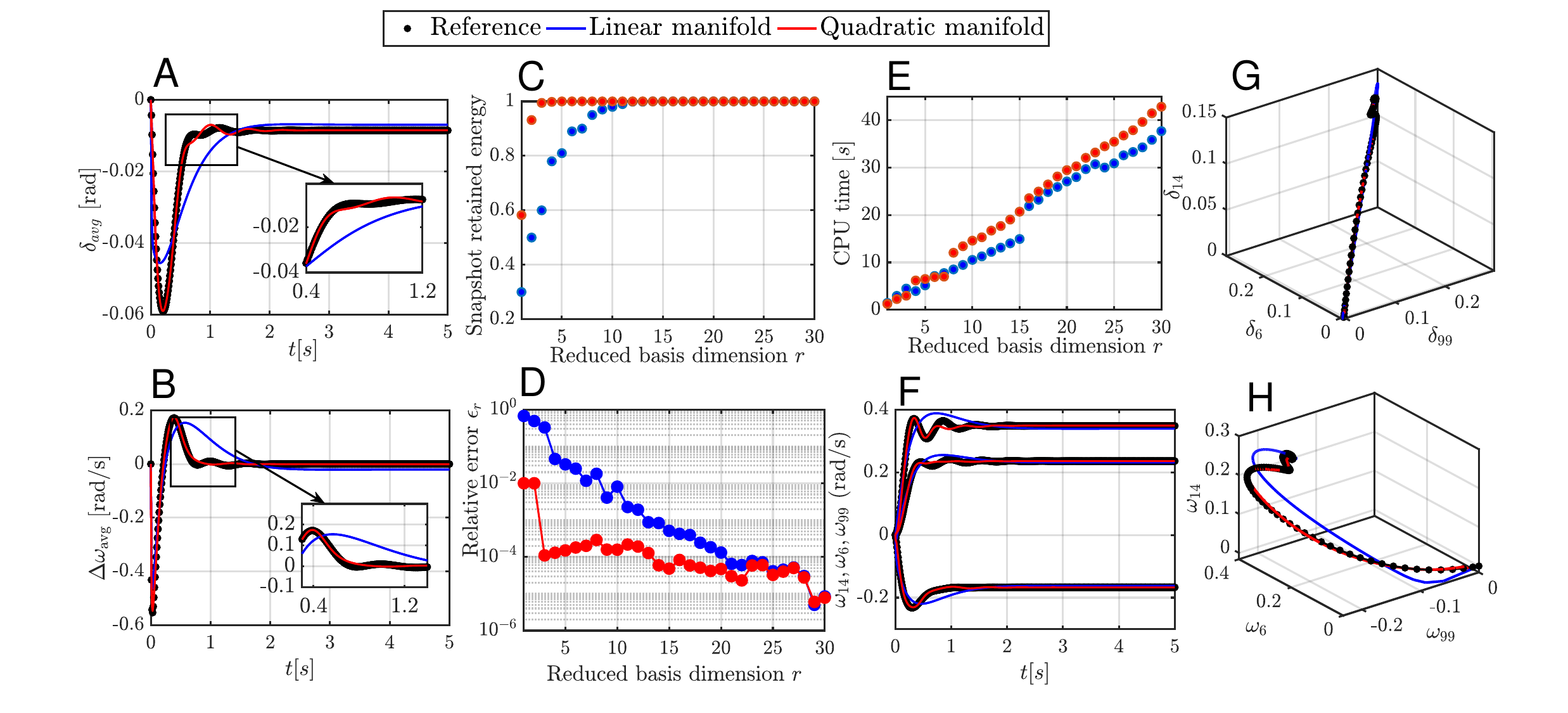} 
	\caption{Simulation results of Polish 2736-bus system: (A) variation of $\delta_{avg}$ with respect to time; (B) variation of $\Delta\omega_{avg}$ with respect to time; (C) snapshot energy captured by the system at various reduced basis dimensions; (D) comparison of relative error for linear and quadratic manifold ROMs; (E) comparison of the CPU times for increasing basis dimension; (F) comparison of the output response for three individual generators with respect to time; (G) comparison os system trajectory onto a linear manifold and quadratic manifold; (H) evolution of angular velocities with respect to time}
	\label{fig:2736}
\end{figure*}
\subsubsection{IEEE 118 Bus System}
The test system comprises 118 buses, 54 generators, and 186 branches. This medium-scale benchmark represents an interconnected transmission network with multiple voltage levels, transformers, and generators, making it suitable for model order reduction studies. A total of $500$ snapshots of rotor angles and frequency deviations were collected over the interval $[0,5s]$. To evaluate the reduction procedures, POD was applied to the snapshot matrix $\mathbf{X}$, while the quadratic manifold was constructed using Algorithm~(\ref{alg:quadratic_manifold}) with the regularization parameter set to $\lambda = 10^{-1}$. The corresponding simulation results are shown in Fig.~(\ref{fig:118}). Plots (A) and (B) illustrate the evolution of the average rotor angle $\delta_{\text{avg}}$ and the average frequency deviation $\Delta\omega_{\text{avg}}$ for both the linear and quadratic manifold techniques. The quadratic manifold closely matches the true system response, whereas the linear manifold exhibits noticeable spurious oscillations. Plot (C) presents the snapshot energy spectra computed using (\ref{eqn-energy_linear}) and (\ref{eqn-energy-nonlinear}) as functions of the reduced basis dimension $r$. We observe that the quadratic manifold retains a higher fraction of the total snapshot energy for a given $r$ compared to the linear subspace. Plot (D) shows the variation of the relative error with respect to basis dimension $r$. We see that the quadratic manifold achieves a more rapid error reduction as $r$ increases, outperforming the linear manifold. This is due to the ill-suitedness of a linear manifold in capturing the nonlinear behaviour of the model. Furthermore, plot (E) shows the online CPU elapsed time; as expected, the quadratic manifold incurs a higher online computational cost due to the Kronecker-product evaluations required by the quadratic ROM. Finally, plots (F)–(H) display the rotor angle trajectories of the three largest generators and the corresponding ROM outputs for a basis dimension of $r=6$. As shown, the quadratic manifold yields significantly improved accuracy compared to the linear manifold in all cases. This is also verified from the $\ell_2$-norm errors and CPU times summarized in Tables~(\ref{tab:L2norm unfaulted}) and (\ref{tab:CPUtimes}), respectively.
\begin{table}[ht]
	\caption{Comparison of $\ell_2$-norm errors for Linear and Quadratic Manifold techniques}
    \centering
	\begin{tabular}{lcccc}
		\toprule
		\multirow{2}{*}{Test Case} & \multicolumn{2}{c}{Linear Manifold} & \multicolumn{2}{c}{Quadratic Manifold} \\ \cline{2-5}
		& Error in $\delta$ & Error in $\omega$ & Error in $\delta$ & Error in $\omega$ \\ \midrule
        \begin{tabular}{l}
            IEEE 118-Bus System\\
             ($n=108, r=6, \lambda=10^{-1}$)
        \end{tabular}
         & $3.4 \times 10^{-2}$  & $2.4 \times 10^{-1}$ & $\mathbf{4.2 \times 10^{-3}}$ & $\mathbf{1.5 \times 10^{-2}}$ \\ \hline
        \begin{tabular}{l}
            IEEE 300-Bus System \\
             ($n=138, r=4, \lambda=1$) 
        \end{tabular}
 & $1.1 \times 10^{-1}$  & $6.2 \times 10^{-1}$ & $\mathbf{2.1 \times 10^{-2}}$ & $\mathbf{5.7 \times 10^{-2}}$ \\ \hline
 \begin{tabular}{l}
            Polish 2736-Bus System \\
             ($n=654, r=3, \lambda=10^{-3}$) 
        \end{tabular}
		& $1.1 \times 10^{-1}$ & $1.1\times 10^{-0}$ & $\mathbf{1.4 \times 10^{-2}}$ & $\mathbf{9.4 \times 10^{-2}}$ \\ \bottomrule
	\end{tabular}
    \label{tab:L2norm unfaulted}
\end{table}
\begin{table}[ht]
\caption{Comparison of CPU times for the unfaulted cases}
    \centering
    \begin{tabular}{lccc}
    \toprule
      Test Case & FOM &Linear ROM & Quadratic ROM  \\ \midrule
         IEEE 118-Bus System &$9.35s$& $1.25s$ & $1.52s$\\
         IEEE 300-Bus System &$15.88s$& $0.61s$& $0.71s$\\
         Polish 2736-Bus System &$3163s$  & $4.47s$      & $4.963s$\\
         \bottomrule
    \end{tabular}
    \label{tab:CPUtimes}
\end{table}
\subsubsection{IEEE 300 bus system}
This test system consists of 300 buses, 69 generators, and 411 branches, resulting in a total state dimension of $n=138$. Its larger size and structural complexity make it suitable for model reduction analysis. As in the previous case, 500 snapshots of the rotor angle and frequency deviations were collected. The numerical simulation results for this case are presented in Fig.~(\ref{fig:300}). Plot (A-B) shows the average system response obtained from the reduced models. The quadratic manifold closely tracks the full-order trajectories, whereas the linear manifold exhibits noticeable deviations. Plot (C) illustrates the retained snapshots energy, computed using (\ref{eqn-energy_linear}) and (\ref{eqn-energy-nonlinear}). For a fixed reduced basis dimension $r$, the quadratic manifold captures a higher fraction of the total snapshot energy than the linear POD subspace for $r\leq 9$. Beyond $r=9$, where the retained energy is close to $100\%$, the least squares problem for estimating $\mathbf{\bar{V}}$ becomes under-determined. Plot (D) shows the relative error as a function of $r$ for both linear and quadratic manifolds. Increasing $r$ reduces the error in both cases; however, the POD-based linear manifold yields the least accurate predictions. The quadratic manifold exhibits faster error decay, enabling greater accuracy with a significantly smaller basis dimension. For instance, a quadratic manifold with $r=3$ achieves approximately the same accuracy as a linear subspace model with $r=8$. The CPU times shown in plot (E) indicate that the computational cost increases with $r$ for both methods, with the POD-based ROM being comparatively faster. Plots (F)-(H) display the rotor angles and frequency deviations of the three largest generators, along with the ROM reconstructions for $r=4$. The quadratic manifold provides a substantially more accurate representation of the rotor dynamics, while the linear manifold exhibits noticeable inaccuracies. These observations are consistent with the $\ell_2$-norm errors in Table (\ref{tab:L2norm unfaulted}) and CPU times reported in Table (\ref{tab:CPUtimes}).
\subsubsection{Polish 2736-bus system}
The model features 2,736 buses, 327 generators, and 3,506 branches, resulting in a total state dimension of $n=654$.  The corresponding simulation results are shown in Fig.~(\ref{fig:2736}). From the figure, plots (A) and (B) display the average rotor angle and frequency deviations, along with the reconstructed responses from the linear and quadratic reduced models for a basis dimension of $r=3$. We observe that the quadratic manifold yields visibly improved agreement with the full-order trajectories. We further see from plot (C) that the total snapshot energy captured by the quadratic manifold is higher than the linear subspace ROM for $r\leq15$. Beyond $r = 15$, the quadratic ROM and the POD ROM capture the same amount of snapshot energy. The relative error trends in plot (D) show that the quadratic ROM error plummets sharply for $r=20$ compared to that of the linear POD ROM. This shows the superior performance of the quadratic manifold in capturing the system's nonlinear behavior more accurately and efficiently. In this case, the quadratic manifold with $r=3$ achieves approximately the same accuracy as a linear subspace model with $r=20$. The flattening of the quadratic manifold error curve for $r>5$ arises from the conditioning of the snapshot matrix $\mathbf{X}$, which introduces numerical errors in less significant POD modes. Plot (E) shows the online CPU times for increasing basis dimensions, and plots (F), (G), and (H) show the rotor responses of the three most significant generators together with the ROM reconstructions. These results further confirm the superior accuracy of the quadratic manifold over the linear subspace ROM. The corresponding relative errors are listed in Table~(\ref{tab:L2norm unfaulted}).  
\subsection{Case II: System with fault behavior}
To assess the robustness of the quadratic ROM under more complex situations, we introduced faults in the IEEE 118-bus EN model. All other parameters--including the reduced basis dimension $r$, regularization parameter $\lambda$, and time-step $\Delta t$, were kept identical to the unfaulted case. Two disturbances were examined: (i) a generator outage and (ii) a three-phase fault followed by line tripping.
\subsubsection{Fault case 1: A generator outage} 
\begin{figure}[htbp]
\centering
\includegraphics[width=0.9\linewidth]{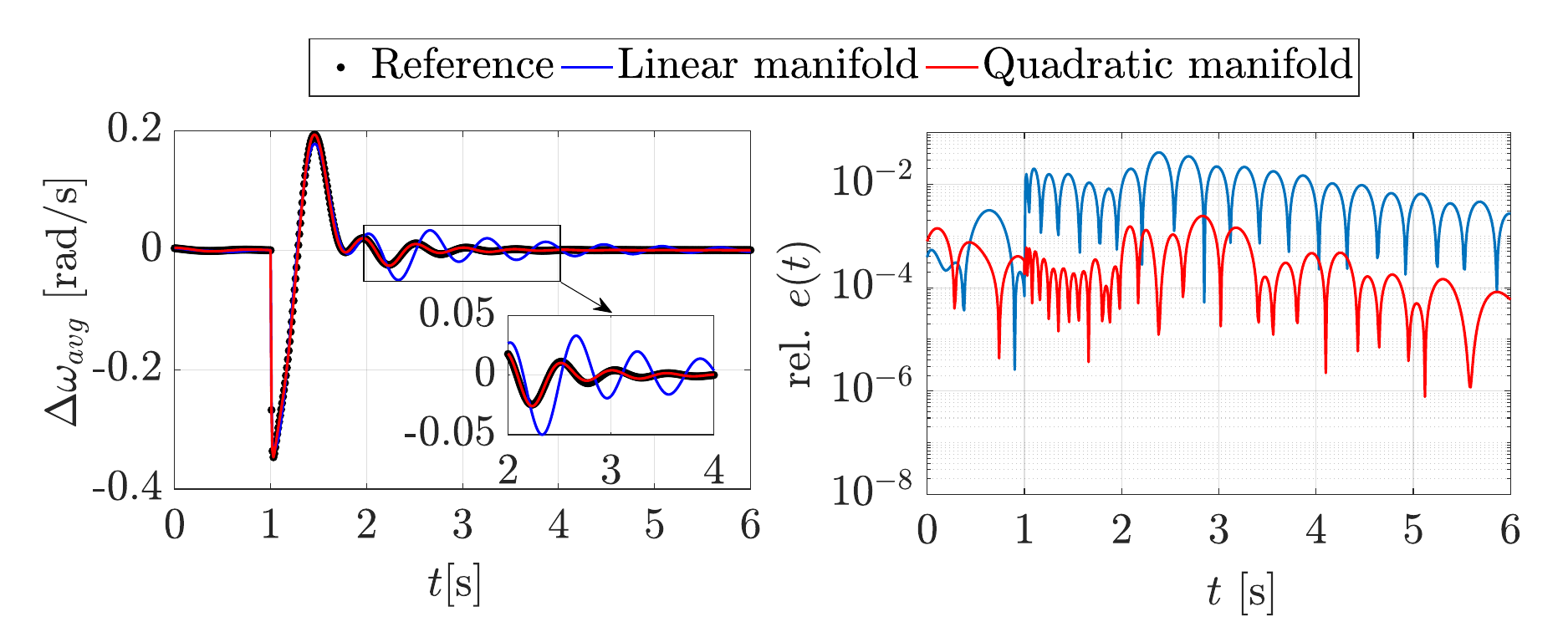}
\caption{Simulation results for fault case 1: (Left) average deviation of $\Delta\omega(t)$; (Right) comparison of relative errors for linear and quadratic manifolds}
\label{fig-fault1}
\end{figure}

In the first scenario, a sudden power outage was simulated by tripping the generator connected to bus 89. The fault was applied at $t=1$s and cleared at $t=1.5$s, resulting in a fault duration of $0.5s$. The corresponding system response is shown in Fig.~(\ref{fig-fault1}). The average frequency exhibits a sharp decline at the fault onset, followed by recovery after clearance.
Both reduced models reproduce the initial transient and the final setting behavior; however, the quadratic manifold ROM provides noticeably closer agreement with the full-order dynamics throughout the post-fault oscillations. This is particularly evident in the inset view over $2$-$4s$ where the quadratic ROM aligns more tightly with the FOM trajectory. Its error remains consistently lower than that of the linear ROM across both transient and steady-state intervals, indicating improved capability to capture the nonlinear behavior during disturbances. These observations are supported by the quantitative error and CPU times reported in Table~(\ref{tab:faultcase 1}). 
\subsubsection{Fault case 2: A three-phase fault with line tripping}
\begin{figure}[!h]
    \centering
~~~~~~\begin{subfigure}[t]{0.67\textwidth}
\includegraphics[width=\linewidth]{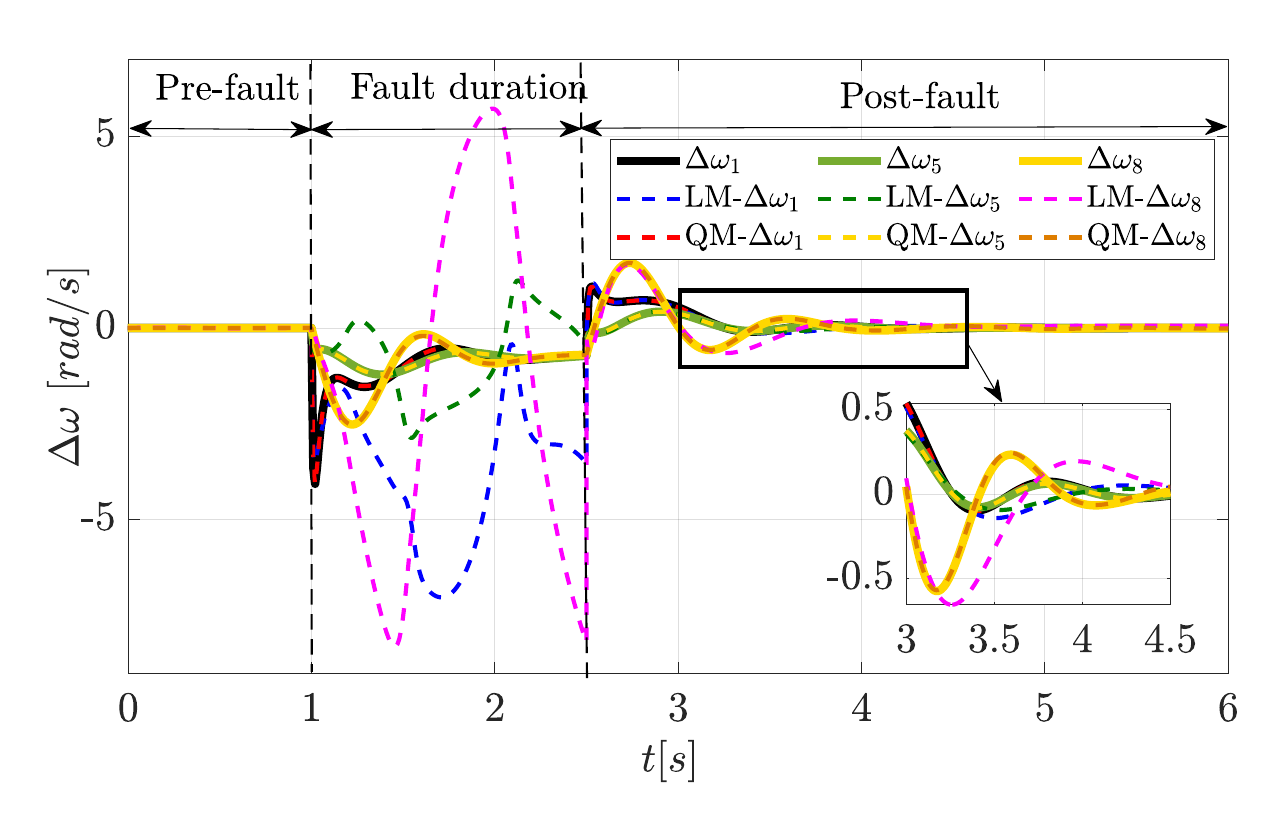}
\end{subfigure}
\begin{subfigure}[t]{0.72\textwidth}
\includegraphics[width=\linewidth]{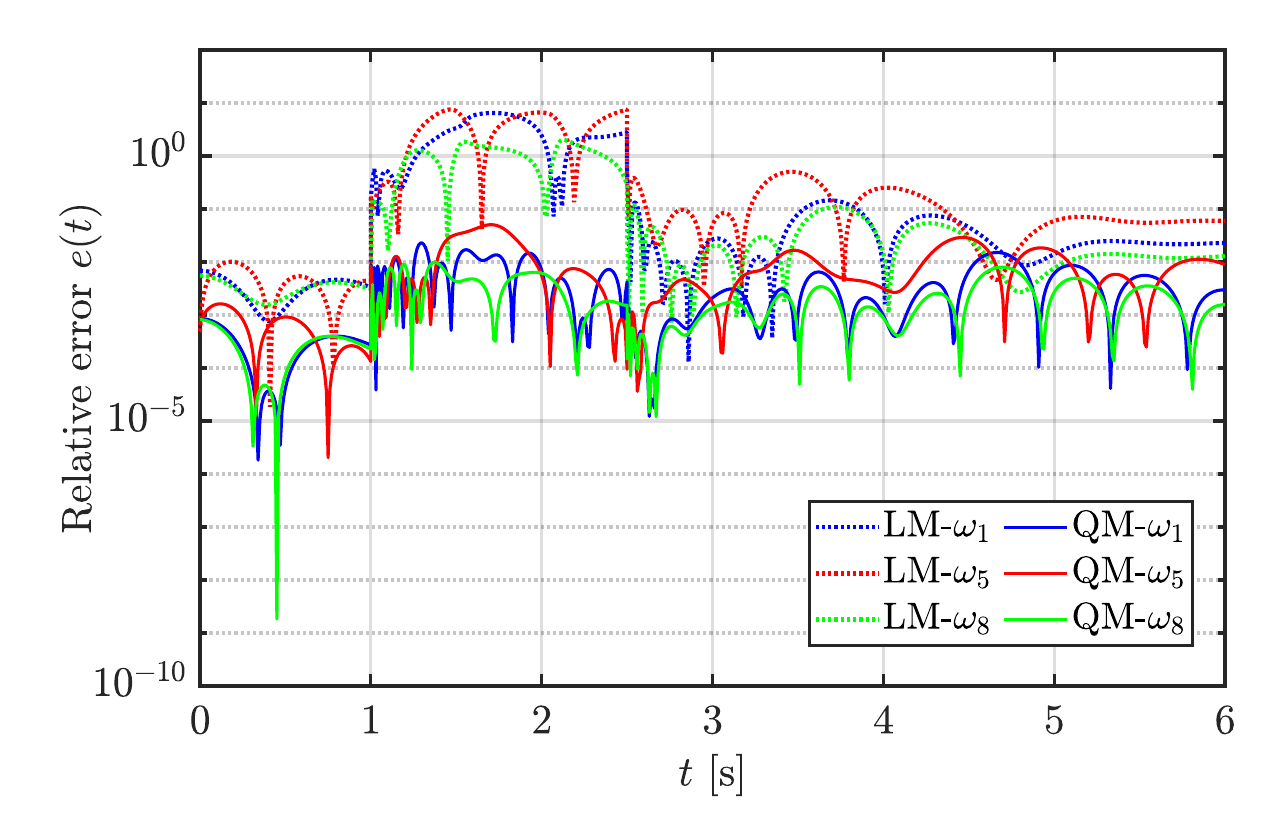}
\end{subfigure}
\caption{Simulation results for  fault case 2: (Top:) comparison of individual frequency deviations $\Delta\omega(t)$ for the $1^{st}$, $5^{th}$ and $8^{th}$ generator of the IEEE 118 Bus system. The solid lines denote the true response, and the dashed lines denote the approximate response from linear and quadratic manifolds. (Bottom:) relative errors $e(t)$ for the $1^{st}$, $5^{th}$, and $8^{th}$ generators of the IEEE 118 Bus system. Solid lines are the error for the quadratic manifold, and dashed lines are the corresponding errors from the linear manifold solution.}
\label{fig-fault2}
\end{figure}

In the second scenario, a three-phase fault was applied at bus 1 at $t=1s$, followed by tripping of the line connecting buses~1 and~2 at $t=2.5s$. The system response and ROM reconstructions are shown in Fig.~(\ref{fig-fault2}). The top plot illustrates the individual frequency deviations of the three generators located closest to the fault. In the top plot, the solid lines represent the FOM response, while the dashed curves indicate the approximations from the linear manifold (LM) and the quadratic manifold (QM) reduced models. All three nearby generators experience a pronounced transient at the fault injection, followed by oscillators' recovery after clearance. The inset highlights the post-fault interval ($3$-$4.5s$), where the quadratic manifold more accurately traces the actual trajectories, capturing both the oscillation envelope and the damping behavior. The linear manifold ROM exhibits larger oscillations and completely fails to capture the strong nonlinear transient behavior caused by the severe fault. The bottom plot reports the relative errors for the same generators. The quadratic manifold consistently achieves lower errors throughout the simulation, while the most significant improvements occur immediately after the disturbance ($1$-$3s$). The corresponding $\ell_2$-norm errors and CPU-time comparisons are provided in Table~(\ref{tab:faultcase 2}). 
\begin{table}[htbp]
	\centering
	\caption{Comparison of CPU times and relative $\ell_2$-norm errors for linear and quadratic manifold techniques for fault case 1}
	\begin{tabular}{lccc}
		\toprule
		{Method} & {CPU Time [$s$]} & {Error in $\delta$} & {Error in $\omega$} \\ \midrule
		FOM & $15.08s$ & -- & -- \\ 
		Linear Manifold  & $1.71s$ & $4.0 \times 10^{-2}$ & $8.1 \times 10^{-2}$ \\ 
		Quadratic Manifold  & $2.12s$ & $\mathbf{4.6 \times 10^{-3}}$ & $\mathbf{4.7 \times 10^{-3}}$ \\ \bottomrule
	\end{tabular}%
    \label{tab:faultcase 1}
	\end{table}
\begin{table}[htbp]
	\centering
	\caption{Comparison of CPU Times and relative $\ell_2$-norms for linear and quadratic manifold techniques for fault case 2}
	\label{tab:cpu_norm}
	\begin{tabular}{lcccc}
		\hline
		\multirow{2}{*}{{Method}} & \multirow{2}{*}{{CPU Time $[s]$}} & \multicolumn{3}{c}{rel. $\ell_2$-norm error} \\
		\cline{3-5}
		& & $\Delta\omega_{1}$ & $\Delta\omega_{5}$ & $\Delta\omega_{8}$ \\
		\midrule
		FOM & $25.58s$  & - & -  & - \\
		Linear Manifold & $1.18s$ & $2.57$ & $3.04$ & $1.22$ \\
		Quadratic Manifold & $2.43s$ & $\mathbf{9.2\times 10^{-3}}$ & $\mathbf{1.9 \times 10^{-2}}$ & $\mathbf{6.5 \times 10^{-3}}$ \\
		\bottomrule
	\end{tabular}%
    \label{tab:faultcase 2}
     \end{table}
\section{Conclusion}\label{sec:con}
This paper presents a nonlinear model reduction technique for large-scale power systems based on the construction of quadratic manifolds. The proposed method learns a nonlinear mapping of the high-dimensional state space by solving a regularized least-squares regression problem using only state measurement data. Unlike classical linear subspace approaches such as Proper Orthogonal Decomposition, the quadratic manifold formulation mitigates the projection errors inherent in linear bases by more accurately capturing the underlying solution manifold. The technique is evaluated on benchmark power system models of increasing size and complexity. The numerical results show that the quadratic manifold consistently outperforms the linear POD method in representing the snapshot energy spectra and reconstructing the system trajectories. Error-norm comparisons further demonstrate that the nonlinear projection yields substantially improved fidelity, particularly under large-disturbance operating conditions, where linear approximations deteriorate. Overall, the proposed method offers a compact, data-driven, and highly accurate reduced-order representation of nonlinear power system dynamics. Its enhanced reconstruction capability makes it a promising tool for transient stability assessment and fast dynamic simulations. Future work may extend the framework to higher-order nonlinear embeddings or real-time implementations for online stability monitoring and control. 

\section{Data Availability} The source code to reproduce the results is available freely at GitHub \url{https://github.com/DanishRaf32/PGmanifold}

\bibliographystyle{unsrt}
\bibliography{references.bib}
\end{document}